\documentstyle[12pt,twoside,leqno]{article}
\newcommand{\eqref}[1]{(\ref{#1})}
\def\myskip{\vspace{8pt}}
\newcommand{\bEq}{\begin{eqnarray}}
\newcommand{\eEq}{\end{eqnarray}}
\newcommand{\beq}{\begin{eqnarray*}}
\newcommand{\eeq}{\end{eqnarray*}}
\newtheorem{Assumption}{\indent Assumption}[assump]
\newtheorem{Definition}{\indent Definition}[section]
\newtheorem{Lemma}{\indent Lemma}[section]
\newtheorem{Proposition}{\indent Proposition}[section]
\newtheorem{Theorem}{\indent Theorem}[section]
\newtheorem{Corollary}{\indent Corollary}[section]
\newtheorem{Remark}{\indent Remark}[section]
\newtheorem{Example}{\indent Example}[section]
\newtheorem{Problem}{\indent Problem}[section]
\newcommand{\bAs}{\begin{Assumption}\em}
\newcommand{\eAs}{\end{Assumption}}
\newcommand{\bDf}{\begin{Definition}\em}
\newcommand{\eDf}{\end{Definition}}
\newcommand{\bLm}{\begin{Lemma}}
\newcommand{\eLm}{\end{Lemma}}
\newcommand{\bPr}{\begin{Proposition}}
\newcommand{\ePr}{\end{Proposition}}
\newcommand{\bTh}{\begin{Theorem}}
\newcommand{\eTh}{\end{Theorem}}
\newcommand{\bCr}{\begin{Corollary}}
\newcommand{\eCr}{\end{Corollary}}
\newcommand{\bRm}{\begin{Remark}\em}
\newcommand{\eRm}{\end{Remark}}
\newcommand{\bEx}{\begin{Example}\em}
\newcommand{\eEx}{\end{Example}}
\newcommand{\bPb}{\begin{Problem}\em}
\newcommand{\ePb}{\end{Problem}}
\newcommand{\bPf}{\par\vspace*{-4pt}\indent{\sc Proof.}\enskip}
\newcommand{\ePf}{\medskip}
\def\QED{\hskip0.1em\hfill\null\ \null\nobreak\hfill\kern3pt\vbox{\hrule\hbox
   {\vrule\kern1pt\vbox{\kern1.7pt\hbox{$\scriptscriptstyle{QED}$}
    \kern0.2pt}\kern1pt\vrule}\hrule}}

\def\END{\hskip0.1em\hfill\null\ \null\nobreak\hfill\kern3pt\vbox{\hrule\hbox
   {\vrule\kern1pt\vbox{\kern1.7pt\hbox{$\,\,\,\vspace{5pt}$}
    \kern0.2pt}\kern1pt\vrule}\hrule}}
\title{\large{{\bf The Prolongation Problem for the Heavenly Equation}}}
\author{{\normalsize M. Palese}
\thanks{Speaker at the Conference. Partially supported by GNFM of CNR, MURST, University of Turin.}
\\{\footnotesize Department of Mathematics,
University of Turin}
\\{\footnotesize Via C. Alberto 10, 10123 Turin, Italy}\\ 
{\footnotesize E--mail: {\sc palese@dm.unito.it}}\\
and\\
{\normalsize R.A. Leo, G. Soliani}
\thanks{Partially supported by INFN, MURST, University of Lecce.}
\\{\footnotesize Department of Physics,
University of Lecce}
\\{\footnotesize Via Arnesano, 73100 Lecce, Italy}\\ 
{\footnotesize E--mail: {\sc leora@le.infn.it, soliani@le.infn.it}}}
\date{}
\begin{document}

\maketitle
\small{
\begin{abstract}
\noindent We provide an exact regular solution of an operator system arising as the prolongation
structure associated with the heavenly equation. This solution is expressed in terms of
operator Bessel coefficients.
\end{abstract}


}

\section{Introduction}

It is well known that the so--called heavenly equation
\bEq\label{1}
u_{xx}+u_{yy}+(e^{u})_{zz}=0\,,
\eEq
where $u=u(x,y,z)$ and subscripts mean partial derivatives, 
occurs in the study of heavenly spaces (Einstein spaces with one rotational Killing vector)
\cite{Boyer,Gegenberg,Gibbons,Lebrun,Plebanski} and in the context of extended conformal 
symmetries \cite{Park}. Reduced versions of Eq.
\eqref{1} have been found in \cite{alfinito}
via the symmetry approach, providing instanton and meron--like
configurations. Furthermore, an algebra of the Virasoro type without central
charge was associated with the heavenly equation by resorting to its invariance under
conformal transformations \cite{alfinito}. 

To investigate the algebraic aspects of the integrability properties of Eq. \eqref{1}, we
shall apply the prolongation technique, which could provide
the relative linear spectral problem
\cite{estabrook,morris,tondo,wahlquist}. In more than two independent variables, the extension of
the prolongation procedure is generally nontrivial and some aspects remain to be explored (see, for
example, \cite{morris,tondo}). However, the main result  of this note is represented by the exact
solution of the prolongation  system associated with Eq. \eqref{1}. This system consists of three
operator  equations in the form of commutator relations which can be written as second order
(operator) differential equations resembling formally equations of the Bessel type. Our solution
is achieved through a series expansion defining operator Bessel coefficients.

\section{The prolongation structure}

Let us introduce on a manifold with local coordinates
$(x,y,z,u,p,q,r)$ the closed differential ideal defined by the set of $3$--forms:
\bEq
\theta_{1} &=& du \wedge dx \wedge dy - rdx \wedge dy \wedge 
dz\,,\label{ideal1}\\
\theta_{2} &=& du \wedge dy \wedge dz - pdx \wedge dy \wedge 
dz\,,\label{ideal2}\\
\theta_{3} &=& du \wedge dx \wedge dz + qdx \wedge dy \wedge 
dz\,,\label{ideal3}\\
\theta_{4} &=& dp \wedge dy \wedge dz - dq \wedge dx \wedge 
dz + \,\label{ideal4}\\
&\hphantom{=}& e^{u}dr\wedge dx \wedge dy + e^{u}r^{2}dx \wedge dy
\wedge dz\,,\nonumber
\eEq 
where $\wedge$ stands for the exterior product.

It is easy to verify the following
\bPr
On every integral submanifold defined by $u = u(x,y,z)$, 
$p=u_{x}$, $q=u_{y}$, $r=u_{z}$, with $dx\wedge dy\wedge dz\neq 0$, the ideal
\eqref{ideal1}--\eqref{ideal4} is equivalent to Eq. \eqref{1}. 
\ePr

Now let us consider the $2$--forms:
\bEq\label{prolongation}
&\hphantom{=}&\Omega^{k} = H^{k}(u, u_{x}, u_{y}, u_{z}; \xi^{m})dx \wedge dy + 
F^{k}(u, u_{x},u_{y}, u_{z}; \xi^{m})dx \wedge dz \\
&\hphantom{=}&+ G^{k}(u, u_{x}, u_{y}, u_{z}; \xi^{m})dy \wedge dz
+ A^{k}_{m}d \xi^{m} \wedge dx+ B^{k}_{m}d \xi^{m} \wedge dz + d \xi^{k}
\wedge dy\,,\nonumber
\eEq 
where $\xi=\{\xi^{m}\}$, $k,\,m\,=\,1,2,\ldots,
{\rm N}$ (N arbitrary), and $H^{k}$, $F^{k}$ and $G^{k}$ are, respectively, the
pseudopotential and functions to be determined. Furthermore, 
the quantities $A^{k}_{m}$
and $B^{k}_{m}$ denote the elements of two $N \times N$ 
constant regular matrices, and
the summation convention over repeated indices is understood.
\bDf
The forms $\Omega^k$ are called the {\em prolongation forms} associated with Eq. \eqref{1}.
Let $\cal{I}$ be the {\em prolonged} ideal generated by the forms 
$\theta_{j}$ and $\Omega^{k}$. 
We say that $\cal{I}$ is closed if $d\Omega^{k}\in \cal{I}$ $(\theta_{j},\Omega^{k})$.
\END\eDf 

\bLm\label{integrability}
The closure condition for $\cal{I}$ yields
\bEq 
&\hphantom{=}& H^{k} = e^{u}u_{z}L^{k}(\xi^{m}) + P^{k}(u,\xi^{m})\,,\label{H}\\
&\hphantom{=}& F^{k} = -\,u_{y}L^{k}(\xi^{m}) + N^{k}(\xi^{m})\,,\label{F}\\
&\hphantom{=}& G^{k} = u_{x}L^{k}(\xi^{m}) + M^{k}(u,\xi^{m})\,,\label{G}
\eEq
where $L^{k}$, $P^{k}$, $N^{k}$, $M^{k}$ are functions of integration. 
\eLm
\bPf
The closure condition is equivalent to the following constraints:
\bEq
&\hphantom{=}& H^{k}_{u_{z}} - e^{u}G^{k}_{u_{x}} = 0\,, \qquad F^{k}_{u_{y}} + G^{k}_{u_{x}} =
0\,,\label{constraints1}\\ 
&\hphantom{=}& H^{k}_{u_{x}} = H^{k}_{u_{y}} = F^{k}_{u_{x}} =  F^{k}_{u_{y}} = G^{k}_{u_{y}} =
G^{k}_{u_{z}} = 0 \,,
\label{constraints4}
\eEq
\bEq
&\hphantom{=}& u_{z}H^{k}_{u} - u_{y}F^{k}_{u} + u_{x}G^{k}_{u} -
 e^{u}u_{z}^{2}G^{k}_{u_{x}} + [G, H]^{k} = 0 \,, \label{structure}\\
&\hphantom{=}& F^{k}_{\xi^{l}} - G^{k}_{\xi^{m}}A^{m}_{l} 
- H^{k}_{\xi^{m}}B^{m}_{l} = 0 \,,\label{spectral1}\\
&\hphantom{=}& F^{k}_{\xi^{m}}(B^{-1})^{m}_{n}G^{n} 
- G^{k}_{\xi^{m}}(B^{-1})^{m}_{n}F^{n} = 0\,,\label{spectral2}
\eEq 
\bEq 
\left[A,B\right] = 0\,,\label{spectral3}
\eEq
\smallskip\\
where $[G, H]^{k} = G^{j}H^{k}_{\xi^{j}} - H^{j}G^{k}_{\xi^{j}}$ (Lie bracket),
$H^{k}_{\xi^{j}} = \frac{\partial H^{k}}{\partial \xi^{j}}$, 
$H^{k}_{u} = \frac{\partial H^{k}}{\partial u}$, and so on.
Equations \eqref{constraints1}--\eqref{constraints4} provide the result.
\QED\ePf

In the following we shall omit the indices $k$, $m$ for simplicity.

\bPr
The following prolongation equations hold
\bEq 
&\hphantom{=}& P_{u} = e^{u}[L,M]\,,\label{P}\\
&\hphantom{=}& M_{u} = - [L,P]\,,\label{M}\\
&\hphantom{=}& \left[M,P\right] = 0\,.\label{comm}
\eEq
\ePr
\bPf
It is a straightforward consequence of Lemma \ref{integrability} (Eq.
\eqref{structure}).
\QED\ePf

\myskip

Hereafter, $L$, $P$, $M$ will be regarded as (regular) operators, in the sense that 
\beq
L \rightarrow L^{j}\frac{\partial}{\partial \xi^{j}}\,, \quad
M \rightarrow M^{j}\frac{\partial}{\partial \xi^{j}}\,, \quad
P \rightarrow P^{j}\frac{\partial}{\partial \xi^{j}}\,,
\eeq
while the Lie brackets become commutators.

\section{Exact solution of the prolongation equations}

Now, let us look for an exact solution to Eqs. \eqref{P}--\eqref{comm}. 

\bDf
For any $A\,=\,A^{j}\frac{\partial}{\partial \xi^{j}}$, we define $\cal{L}[A]$ by
\bEq
{\cal{L}}[A] = [L, A]\,.
\eEq
\END\eDf
 
\bDf
We define the $n$--th power of the operator ${\cal{L}}$ by setting
\bEq
{\cal{L}}^{n} [A] = [L,[L,\dots ,[L, A]\dots]\,,\label{binomial}
\eEq
where $L$ appears $n$--times, and ${\cal L}^{0}[A] = A$.
\END\eDf

\bRm
The prolongation Eqs. \eqref{P} and \eqref{M} can be written as

the second order operator equations 
\bEq 
&\hphantom{=}& P_{tt} - \frac{1}{t}P_{t} + {\cal{L}}^{2} [P] = 0\,,
\label{P2}\\
&\hphantom{=}& M_{tt} + \frac{1}{t}M_{t} + \frac{1}{t}{\cal{L}}^{2} [M] = 0\,,
\label{M2}
\eEq
where $t = 2e^{\frac{u}{2}}$.
\END\eRm

\bRm
Notice that the above equations resemble formally conventional Bessel equations of the type
\beq 
 \kappa_{tt} - \frac{1}{t}\kappa_{t} + \omega^{2} \kappa = 0\,,
\qquad  \chi_{tt} + \frac{1}{t}\chi_{t} + \frac{1}{t}\omega^{2} \chi = 0\,,
\eeq
with $\omega$ a constant, whose regular solutions at $t=0$ are given by
\beq
\kappa(t)=p_{0}\frac{t}{2}J_{1}(t\omega)\,,
 \qquad \chi(t)=m_{0}J_{0}(t\omega)\,,
\eeq
with $p_{0}$, $m_{0}$ constants of integration.
\END\eRm

\myskip

We are interested in solutions of Eqs. \eqref{P2} and \eqref{M2}
which are regular at $t = 0$. We proceed in a heuristic way generalizing the
scheme working out for the solution of the operator equations:  

\bEq
&\hphantom{=}& {\tilde P}_{tt} + {\cal{L}}^{2} [{\tilde P}] = 0\,,\quad {\tilde M}_{tt} +
{\cal{L}}^{2} [{\tilde M}] = 0\,, \label{harmonic2}
\eEq
whose formal solution is
\bEq 
&\hphantom{=}& {\tilde P}(t,\xi) = e^{it{\cal L}}[A_{0}(\xi)] + 
e^{-it{\cal L}}[B_{0}(\xi)]\,,\label{x}\\
&\hphantom{=}& {\tilde M}(t,\xi) = e^{it{\cal L}}[C_{0}(\xi)] +
e^{-it{\cal L}}[D_{0}(\xi)]\,,\label{xx}
\eEq
where $A_0$, $B_0$, $C_0$, $D_0$ can be determined from the initial conditions.

\bDf
Let ${\bf J_0} (\cdot)$ and ${\bf J_1} (\cdot)$ be formally defined by
\bEq
&\hphantom{=}& {\bf J_0} (t{\cal{L}}[M_0]) = 
\sum_{m=0}^{\infty}\frac{(-1)^m}{(m!)^2}
\left(\frac{t}{2}\right)^{2m}{\cal{L}}^{2m}
[M_0]\,,\label{J0cal}\\
&\hphantom{=}& {\bf J_1} (t{\cal{L}} [P_0]) =
\sum_{m=0}^{\infty}\frac{(-1)^m}{m!(m+1)!}
\left(\frac{t}{2}\right)^{1+2m}{\cal{L}}^{1+2m}
[P_0]\,.\label{J1cal}
\eEq
\END\eDf

\bRm
Here ${\bf J_0}$ and ${\bf J_1}$ 
are borrowed by the series expansion of 
the Bessel functions $J_{\nu}$, with $\nu$ an integer:
\bEq 
J_{\nu} (\tau) = \sum_{m=0}^{\infty}\frac{(-1)^m}{m!(m+\nu)!}
\left(\frac{\tau}{2}\right)^{\nu+2m}{\cal{L}}^{\nu+2m}[P_0]\,,
\eEq
for $\nu = 0$ and $\nu = 1$, respectively (see \cite{watson}).
\END\eRm

The following result holds, whose proof is straightforward.
\bPr
A solution of equations \eqref{P2} and \eqref{M2} regular at $t=0$ is given by
\bEq 
P = \frac{t}{2}{\bf J_1} (t{\cal{L}} [P_0])\,,\qquad
M = {\bf J_0} (t{\cal{L}} [M_0])\,,\label{first solution}
\eEq 
where $M_0 \equiv M_0(\xi) = M(t;\xi)\mid _{t=0}$ and
$P_0 \equiv P_0(\xi)$ is such that $[L,P_0] = [L,M_0]$.
\ePr

\bRm
Inserting the operators \eqref{first solution} into 
Eqs. \eqref{P2} and \eqref{M2} yields the condition $[{\cal L}[M_0],M_0] = 0$.
\END\eRm

\bRm
We point out that the search of solutions to Eqs. \eqref{P2} and \eqref{M2}, 
regular at $t=0$, is equivalent to
tackle the corresponding Cauchy problem with the initial conditions
$P(0;\xi) = 0$, $P_{t}(t,\xi)\mid _{t=0} = 0$, $M(0;\xi) = M_0$, $M_{t}(t,\xi)\mid _{t=0}
= 0$.\END
\eRm

\myskip

For pratical pourposes, {\em e.g.} to determine the spectral problem
and B\"acklund transformations associated with Eq. \eqref{1}, we shall rewrite 
the solutions \eqref{first solution} in a form which 
contains the operator $L$ and not the operator ${\cal L}$.
We shall give a regular solution of Eqs. \eqref{P2} and \eqref{M2} following a scheme similar to
that working out for the case \eqref{harmonic2}.

By induction it is easy to prove the following Lemma.

\bLm
For $n>0$, \eqref{binomial} takes the form 
\bEq
{\cal{L}}^{n} [A] = \sum_{k=0}^{n}(-1)^k{n \choose 
k}L^{n-k}AL^{k}\,.\label{fundamental tool}
\eEq
\eLm

As a straightforward consequence we have
\bEq
&\hphantom{=}&e^{it{\cal L}}[A_{0}] =\sum_{n=0}^{\infty}\frac{(it)^{n}}{n!}
\sum_{k=0}^{n}(-1)^k{n \choose k}L^{n-k}[A_0]L^{k}\nonumber\\
&\hphantom{=}&= \sum_{j=0}^{\infty}\frac{(it)^{j}}{j!}L^{j}[A_0]
\left(\sum_{K=0}^{\infty}\frac{(it)^k}{k!}L^{k}\right)\equiv e^{itL}[A_{0}]e^{-itL}\,,
\label{baker}
\eEq
which is just the {\em Baker--Campbell--Hausdorff expansion} (see 
{\em e.g.} \cite{baker}).

Thus, we can rewrite the solutions of the operator equations \eqref{harmonic2}
in a more suitable form.

\subsection{Operator Bessel coefficients}

In order to express the solutions \eqref{first solution} of Eqs. 
\eqref{P2} and \eqref{M2} in terms of the operator $L$, we shall introduce operator
Bessel coefficients by means of a formal expansion analogous to that used in the
case of Bessel functions.

\bDf
Let $X$ be a regular operator.
We define the {\em operator Bessel coefficients} ${\bf J_m}(tX)$, 
as the coefficients of the formal expansion:
\bEq
e^{\frac{t}{2}X(z-1/z)} = \sum_{m=-\infty}^{\infty}z^{m}{\bf J_{m}}(tX)\,.
\label{fundamental expansion}\END
\eEq
\eDf

\bRm\label{laurent}
We stress that the Laurent series on the right side is uniformly convergent.
\END\eRm

First we prove a technical Lemma.
\bLm\label{bessel}
Operator Bessel coefficients satisfy the following recurrence and derivation formulae:
\bEq
&\hphantom{=}& {\bf J_{- k}}(tL) 
= (-1)^{k}{\bf J_{k}}(tL)\,,\label{negative index}\\
&\hphantom{=}& 2k{\bf J_{k}}(tL) = tL[{\bf J_{k-1}}(tL) + {\bf 
J_{k+1}}(tL)]\,,\label{positive recurrence}\\
&\hphantom{=}& 2\frac{d}{dt}[{\bf J_{k}}(tL)]
= L[{\bf J_{k-1}}(tL) - {\bf J_{k+1}}(tL)]\,.\label{negative 
recurrence}\\
&\hphantom{=}& \frac{d}{dt}[t^{k}{\bf J_{k}}(tL)] 
= Lt^{k}{\bf J_{k-1}}(tL)\,,\label{positive derivative}\\
&\hphantom{=}& \frac{d}{dt}[t^{- k}{\bf J_{k}}(tL)] 
= - Lt^{- k}{\bf J_{k+1}}(tL)\,,\label{negative derivative}
\eEq
\eLm

\bPf
By differentiating with respect to $z$ the formal expansion \eqref{fundamental
expansion}, with $X=L$, we obtain
\bEq
\frac{1}{2}tL(1+\frac{1}{z^{2}})
\sum_{k=-\infty}^{\infty}z^{k}{\bf J_{k}}(tL) = 
\sum_{k=-\infty}^{\infty}kz^{k-1}{\bf J_{k}}(tL)\,.\nonumber
\eEq
Then, by equating coefficients of 
$z^{k-1}$ in the above identity we obtain formula \eqref{positive recurrence}.
Furthermore, if we differentiate the formal expansion 
with respect to $t$ we have
\bEq
\frac{1}{2}L(z-1/z)\sum_{k=-\infty}^{\infty}z^{k}{\bf J_{k}}(tL) = 
\sum_{k=-\infty}^{\infty}z^{k}\frac{d}{dt}{\bf J_{k}}(tL)\,.\nonumber 
\eEq 
By equating coefficients of $z^{k}$ on either side of 
this identity we obtain formula \eqref{negative recurrence}. 
Formulae \eqref{positive derivative} and \eqref{negative derivative} can be determined by
adding and substracting 
\eqref{positive recurrence} and \eqref{negative recurrence}, while Eq. \eqref{negative index}
follows  directly from \eqref{fundamental
expansion} and from the hypothesis that $L$ 
is a regular operator.
\QED\ePf

\subsection{A form of the solution of the prolongation equations in terms of $L$}

In the following we shall provide an equivalent
solution to Eqs. \eqref{P2} and \eqref{M2} in terms of $L$ which is 
in some sense an analogue of formula \eqref{baker}.

\bPr
We can rewrite the solution \eqref{first solution} in terms of $L$ as follows:
\bEq
P = \frac{t}{2}\sum_{k=-\infty}^{\infty}{\bf J_{k+1}}(tL)P_{0}{\bf 
J_{k}}(tL)\,, \qquad M = \sum_{k=-\infty}^{\infty}{\bf J_{k}}(tL)M_{0}{\bf 
J_{k}}(tL)\,.\label{M second solution}
\eEq
\ePr

\bPf
To verify that the operators \eqref{M second solution} obey Eqs.
\eqref{P2} and \eqref{M2}, we refer to Lemma \ref{bessel}. 
In fact, by resorting to \eqref{positive derivative}, \eqref{negative derivative} and
\eqref{negative index} we have
\beq
P_{tt} =\frac{1}{2}\sum_{k=-\infty}^{\infty}[{\bf  J_{k}}(tL){\cal L}[P_{0}]{\bf
J_{k}}(tL)]+ \frac{t}{2}\sum_{k=-\infty}^{\infty}
[{\bf J_{k}}(tL){\cal L}^{2}[P_{0}]{\bf J_{k+1}}(tL)]\,.
\eeq
On the other hand, since
\bEq
&\hphantom{=}&{\cal L}^{2}[P] = - \frac{t}{2}\sum_{k=-\infty}^{\infty}(-1)^{k}
[{\bf J_{k}}(tL)L^{2}P_{0}{\bf J_{1-k}}(tL) - 2{\bf J_{k}}(tL)LP_{0}L{\bf J_{1-k}}(tL)\nonumber\\
&\hphantom{=}&+ {\bf J_{k}}(tL)P_{0}L^{2}{\bf J_{1-k}}(tL)] = -
\frac{t}{2}\sum_{k=-\infty}^{\infty} [{\bf J_{k}}(tL){\cal L}^{2}[P_{0}]{\bf
J_{k+1}}(tL)]\,,\nonumber
\eEq
we obtain the result.

In a similar way, by virtue of \eqref{positive recurrence}, we can write
\bEq
L{\bf J_{k+2}}(tL) = \frac{2(k+1)}{t}{\bf J_{k+1}}(tL) - L{\bf 
J_{k}}(tL)\,.
\eEq
This, with the help of Lemma \ref{bessel}, gives the following expression for $M_{tt}$:
\bEq
&\hphantom{=}&-\sum_{k=-\infty}^{\infty}
[{\bf J_{k}}(tL){\cal L}^{2}[M_{0}]{\bf J_{k}}(tL)]+ \sum_{k=-\infty}^{\infty}(-1)^{k}2(k-1)[{\bf
J_{k-1}}(tL)LM_{0}{\bf  J_{-k}}(tL) \nonumber\\
&\hphantom{=}&- 2{\bf J_{k-1}}(tL)M_{0}L{\bf J_{-k}}(tL) + 
{\bf J_{-k}}(tL)M_{0}L{\bf J_{k-1}}(tL)]\,.\nonumber
\eEq
Then, using the expression
\bEq
{\cal L}^{2}[M] = \sum_{k=-\infty}^{\infty}
[{\bf J_{k}}(tL){\cal L}^{2}[M_{0}]{\bf J_{-k}}(tL)]\,,\nonumber
\eEq
and taking into account \eqref{negative index} and \eqref{positive 
recurrence}, we have
\beq
&\hphantom{=}&tM_{tt}+M_{t} +{\cal L}^{2}[M]=\sum_{k=-\infty}^{\infty} [2k{\bf
J_{k-1}}(tL)LM_{0}{\bf J_{k}}(tL)+\\
&\hphantom{=}&- 2(k-1){\bf J_{k-1}}(tL)M_{0}L{\bf J_{k}}(tL)] = 0\,.
\eeq
Hence Eq. \eqref{M2} is satisfied. This completes the proof.\QED
\ePf

\section{Conclusions}
We have solved the prolongation problem for Eq. 
\eqref{1} in terms of a series expansion of
operators which can be interpreted as generalized 
Bessel coefficients. The operators
\eqref{M second solution} 
have been derived under the hypothesis that
$P$ and $M$ are regular at $t=0$, {\em i.e.} at $u\to -\infty$. 
The ``key'' for our result is based
on Eqs. \eqref{P2} and \eqref{M2}, which constitute an extended form 
(of operator Bessel type) of the
operator equations \eqref{harmonic2}. 
We point out that the introduction
of the operator \eqref{fundamental tool} has strongly facilited our task. 
To this regard, a useful
step is represented by Eq. \eqref{baker}, where a 
correspondence is established between the
operators ${\cal L}$ and $L$ via the Baker--Campbell--Hausdorff formula.

In theory, the knowledge  of $H$, $F$ and $G$ (see \eqref{H}--\eqref{G}) 
may be exploited to find
the spectral problem related to Eq. \eqref{1}. However, in opposition 
to what happens in other
cases, the determination of the spectral problem of Eq. 
\eqref{1} within the prolongation
scheme offers notable difficulties, mainly owing to the fact that 
we have been able to get only a
solution of the prolongation equations \eqref{P}--\eqref{comm} which 
is regular at $t=0$. At present the construction of 
a prolongation algebra 
whose elements depend uniquely on the
pseudopotential variables remain an open
problem.

\subsection{Acknowledgments}
One of us (M. P.) would thank M. Ferraris and M. Francaviglia for useful discussions. 
 

\end{document}